\documentclass[11pt]{amsart}

\usepackage{amsfonts,amsmath,amssymb,amsthm}

\newcommand{\Cc}{\mathbb{C}}  

\newcommand{\Nn}{\mathbb{N}}

\newcommand{\defi}[1]{\emph{#1}}

\renewcommand{\epsilon}{\varepsilon}
\renewcommand{\le}{\leqslant}
\renewcommand{\ge}{\geqslant}

\newcommand{\Baff}{{\mathcal{B}_\mathit{\!aff}}}
\newcommand{\Binf}{{\mathcal{B}_\infty}}
\newcommand{\Binfx}{{\mathcal{B}_{\infty,x}}}
\newcommand{\Binfy}{{\mathcal{B}_{\infty,y}}}
\newcommand{\B}{{\mathcal{B}}}

\newcommand{\Aut}{\mathop{\mathrm{Aut}}\nolimits}

\newcommand{\supp}{\mathop{\mathrm{supp}}\nolimits}

{\theoremstyle{plain}
\newtheorem{theorem}{Theorem}    
\newtheorem{lemma}[theorem]{Lemma}        
\newtheorem{corollary}[theorem]{Corollary}        
}
{\theoremstyle{remark}
\newtheorem*{remark*}{Remark}  

}

\begin{document}

\title{Newton polygons and families of polynomials}
\author{Arnaud Bodin}
\address{Laboratoire Paul-Painlev\'e, Math\' ematiques,
Universit\'e de Lille 1, \  59655 Villeneuve d'Ascq, France.}
\email{Arnaud.Bodin@math.univ-lille1.fr}

\begin{abstract}
We consider a continuous family $(f_s)$, $s\in[0,1]$ of complex polynomials in two variables 
with isolated singularities, that are Newton non-degenerate. We suppose that
the Euler characteristic of a generic fiber is constant.
 We firstly prove that the set of critical 
values at infinity depends continuously on $s$,
and secondly that the degree of the $f_s$ is constant 
(up to an algebraic automorphism of $\Cc^2$).
\end{abstract}

\maketitle

\section{Introduction}

We consider a family $(f_s)_{s\in[0,1]}$ of complex polynomials in two variables 
with isolated singularities. We suppose that coefficients are continuous functions
of $s$. For all $s$, there exists a finite \defi{bifurcation set} $\B(s)$ such that 
the restriction of $f_s$ above $\Cc\setminus \B(s)$ is a locally trivial fibration.
It is known that $\B(s) = \Baff(s) \cup \Binf(s)$, where $\Baff(s)$ is the set of
\defi{affine critical values}, that is to say the image by $f_s$ of the critical points;
$\Binf(s)$ is the set of \defi{critical values at infinity}.
For $c \notin \B(s)$, the Euler characteristic verifies 
$\chi(f_s^{-1}(c)) = \mu(s) + \lambda(s)$,
where $\mu(s)$ is the \defi{affine Milnor number} and $\lambda(s)$ is the \defi{Milnor number
at infinity}.

We will be interested in families such that the sum  $\mu(s)+\lambda(s)$ is constant.
These families are interesting in the view of $\mu$-constant type theorem, see \cite{HZ,HP,Ti,Bo,BT}.
We say that a multi-valued function $s\mapsto F(s)$ is \defi{continuous} if
at each point $\sigma \in [0,1]$ and at each value $c(\sigma)\in F(\sigma)$
 there is a neighborhood $I$ of $\sigma$ such that
for all $s \in I$, there exists $c(s)\in F(s)$ near $c(\sigma)$.
$F$ is \defi{closed}, 
if, for all points $\sigma \in [0,1]$, 
for all sequences $c(s)\in F(s)$, $s\not=\sigma$, such that $c(s) \rightarrow c(\sigma) \in \Cc$
as $s\rightarrow \sigma$, then $c(\sigma) \in F(\sigma)$.
It it is well-known that $s \mapsto \Baff(s)$ is a continuous multi-valued function.
But it is not necessarily closed: for example $f_s(x,y) = (x-s)(xy-1)$, then
for $s\not=0$, $\Baff(s)=\{ 0, s \}$ but $\Baff(0)=\varnothing$.

We will prove that $s \mapsto \Binf(s)$ and $s \mapsto \B(s)$ are closed continuous functions
 under some assumptions.

\begin{theorem}
\label{th:cont}
Let $(f_s)_{s\in[0,1]}$ be a family of complex polynomials
such that $\mu(s)+\lambda(s)$ is constant and such that 
$f_s$ is (Newton) non-degenerate for all $s\in[0,1]$, then
the multi-valued function $s \mapsto \Binf(s)$ is continuous 
and closed.
\end{theorem}

\begin{remark*}
As a corollary we get the answer to a question of D.~Siersma: is it possible to
find a family with $\mu(s)+\lambda(s)$ constant such that $\lambda(0) > 0$
(equivalently $\Binf(0)\neq \varnothing$) and $\lambda(s)=0$ (equivalently $\Binf=\varnothing$)
for $s\in ]0,1]$? Theorem \ref{th:cont} says that it is not possible for non-degenerate polynomials.
Moreover for a family  with $\mu(s)+\lambda(s)$ constant and $\lambda(s)>0$ for $s\in]0,1]$ 
we have $\lambda(0) \ge \lambda(s) > 0$ by the (lower) semi-continuity of $\mu(s)$.
In the case of a $\mathcal{F}ISI$ deformation of polynomials of constant degree with a non-singular total
space, the answer can be deduced from \cite[Theorem 5.4]{ST}.
\end{remark*}

\begin{remark*}
Theorem \ref{th:cont} does not imply that $\mu(s)$ and $\lambda(s)$ are constant.
For example let the family $ f_s(x,y) = x^2y^2+sxy+x.$
Then  for $s=0$, $\mu(0)=0$, $\lambda(0) = 2$ with $\Binf(0)=\{0\}$, and for $s\not=0$,
$\mu(s)=1$, $\lambda(s)=1$ with $\Baff(s)=\{0\}$ and $\Binf(s) = \{ -\frac{s^2}{4} \}$.
\end{remark*}

The multi-valued function $s\mapsto \Baff(s)$ is continuous but not necessarily 
 closed even if $\mu(s)+\lambda(s)$ is constant, for example (see \cite{Ti}):
$f_s(x,y)=x^4-x^2y^2+2xy+sx^2$, then $\mu(s)+\lambda(s) = 5$.
We have $\Baff(0)=\{0\}$, $\Binf(0)=\{1\}$ and for $s\not=0$, $\Baff=\{0,1-\frac{s^2}{4}\}$, $\Binf(s)= \{1\}$.
We notice that even if $s\mapsto \Baff(s)$ is not closed, the map $s\mapsto \B(s)$ is closed.
This is expressed in the following corollary (of Theorems \ref{th:cont} and \ref{th:deg}):
\begin{corollary}
\label{th:contb}
Let $(f_s)_{s\in[0,1]}$ be a family of complex polynomials
such that $\mu(s)+\lambda(s)$ is constant and such that 
$f_s$ is non-degenerate for all $s\in[0,1]$.
Then the multi-valued function $s \mapsto \B(s)$ is continuous 
and closed.
\end{corollary}

We are now interested in the constancy of the degree; in all hypotheses of global 
$\mu$-constant theorems the degree of the $f_s$ is supposed not to change (see \cite{HZ,HP,Bo,BT})
and it is the only non-topological hypothesis.  We prove that for non-degenerate 
polynomials in two variables the degree is constant except for a few cases, where
the family is of quasi-constant degree.
We will define in a combinatoric way in paragraph \ref{sec:demdeg} what a family of \defi{quasi-constant degree}
is, but the main point is to know that such a family is of constant degree
 up to some algebraic automorphism of $\Cc^2$. More precisely, for each value $\sigma\in[0,1]$
there exists $\Phi \in \Aut \Cc^2$ such $f_s \circ \Phi$ is of constant degree, for $s$ 
in a neighborhood of $\sigma$. For example 
the family $f_s(x,y)=x+sy^2$ is of quasi-constant degree while 
the family $f_s(x,y)=sxy+x$ is not.

\begin{theorem}
\label{th:deg}
Let $(f_s)_{s\in[0,1]}$ be a family of complex polynomials
such that $\mu(s)+\lambda(s)$ is constant and such that 
$f_s$ is non-degenerate for all $s\in]0,1]$, then
either $(f_s)_{s\in[0,1]}$ is of constant degree or $(f_s)_{s\in[0,1]}$ 
is of quasi-constant degree.
\end{theorem}

\begin{remark*}
In theorem \ref{th:deg}, $f_0$ may be degenerate.
\end{remark*}

As a corollary we get a $\mu$-constant theorem without hypothesis on the degree:
\begin{theorem}
\label{th:mucst}
Let $(f_s)_s\in[0,1]$ be a family of polynomials in two variables with 
isolated singularities such that the coefficients are continuous function of $s$.
We suppose that $f_s$ is non-degenerate for $s\in ]0,1]$, and that 
the integers $\mu(s)+\lambda(s), \# \B(s)$ are constant ($s\in [0,1]$)
then the polynomials $f_0$ and $f_1$ are topologically equivalent.
\end{theorem}
It is just the application of the $\mu$-constant theorem of \cite{Bo}, \cite{BT} 
to the family $(f_s)$ or $(f_s\circ \Phi)$.
Two kinds of questions can be asked : are Theorems \ref{th:cont} and \ref{th:deg} 
true for degenerate polynomials? are they true for polynomials in more than $3$
variables?
I would like to thank Prof.~G\"{u}nter Ewald for
discussions concerning Theorem \ref{th:deg} in $n$ variables
(that unfortunately only yield that the given proof cannot be easily generalized).

\section{Tools}

\subsection{Definitions}

We will recall some basic facts about Newton polygons, see \cite{K},
\cite{CN}, \cite{NZ}.
Let $f \in \Cc[x,y]$, $f(x,y) = \sum_{(p,q)\in \Nn^2} a_{p,q}x^py^q$. We denote 
$\supp(f) = \{ (p,q) \mid a_{p,q} \not=0 \}$, by abuse $\supp(f)$ will also denote
the set of monomials $\{ x^py^q \mid (p,q)\in\supp(f)\}$.
$\Gamma_-(f)$ is the convex closure
of $\{(0,0)\}\cup \supp(f)$, $\Gamma(f)$  is the union of closed faces which do not contain
$(0,0)$. 
For a face $\gamma$, $f_\gamma = \sum_{(p,q) \in \gamma} a_{p,q}x^py^q$. The polynomial
$f$ is \defi{(Newton) non-degenerate} if for all faces $\gamma$ of $\Gamma(f)$ the system
$$\frac{\partial f_\gamma}{\partial x}(x,y) = 0; \quad \frac{\partial f_\gamma}{\partial y}(x,y) = 0$$
has no solution in $\Cc^*\times \Cc^*$.

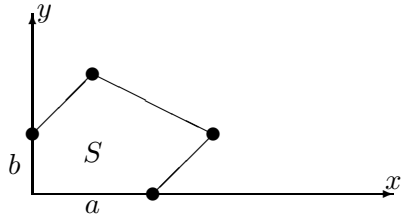
\begin{figure}[ht]
\begin{center}
\unitlength 0.8mm
\begin{picture}(60,30)(0,-3)

\put(0,0){\vector(1,0){60}}
\put(0,0){\vector(0,1){30}}
\put(60,2){\makebox(0,0){$x$}}
\put(2,30){\makebox(0,0){$y$}}

\put(0,10){\circle*{2}}
\put(0,10){\line(1,1){10}}
\put(10,20){\circle*{2}}
\put(10,20){\line(2,-1){20}}
\put(30,10){\circle*{2}}
\put(30,10){\line(-1,-1){10}}
\put(20,0){\circle*{2}}

\put(10,7){\makebox(0,0){$S$}}
\put(10,-2){\makebox(0,0){$a$}}
\put(-3,5){\makebox(0,0){$b$}}
\end{picture}

\caption{Newton polygon of $f$ and $\nu(f) = 2S-a-b+1$. \label{fig:nu}}
\end{center}
\end{figure}

We denote by $S$ the area of $\Gamma_-(f)$, by $a$ the length of the intersection of $\Gamma_-(f)$
with the $x$-axis, and by $b$ the length of the intersection of $\Gamma_-(f)$ with the $y$-axis 
(see Figure \ref{fig:nu}).
We define 
$$\nu(f) = 2S-a-b+1.$$

\subsection{Milnor numbers}

The following result is due to Pi.~Cassou-Nogu\`es \cite{CN}, it is an improvement of
Kouchnirenko's result.

\begin{theorem}
\label{th:CN}
Let $f \in \Cc[x,y]$ with isolated singularities. Then
\begin{enumerate}
  \item $\mu(f) + \lambda(f) \le \nu(f)$.
  \item If $f$ is non-degenerate then $\mu(f) + \lambda(f) = \nu(f)$.
\end{enumerate}
\end{theorem}

\subsection{Critical values at infinity}

We recall the result of A.~N\'em\'ethi and A.~Zaharia on
how to estimate critical values at infinity.
A polynomial $f\in \Cc[x,y]$ is \defi{convenient for the $x$-axis} if there exists a monomial
$x^a$ in $\supp(f)$  ($a>0$); $f$ is \defi{convenient for the $y$-axis} 
if there exists a monomial $y^b$ in $\supp(f)$ ($b>0$) ; $f$ is \defi{convenient}
if it is convenient for the $x$-axis and the $y$-axis.
It is well-known (see \cite{Br}) that:
\begin{lemma}
A non-degenerate and convenient polynomial with isolated singularities
has no critical value at infinity: $\Binf = \varnothing$.
\end{lemma}

Let $f\in \Cc[x,y]$ be a  polynomial with $f(0,0)=0$ not depending only on one variable.
Let $\gamma_x$ and $\gamma_y$ the two faces of $\Gamma_-(f)$ that contain the origin.
If $f$ is convenient for the $x$-axis then we set $\mathcal{C}_x = \varnothing$
otherwise $\gamma_x$ is not included in the $x$-axis and we set 
$$\mathcal{C}_x = \bigg\lbrace f_{\gamma_x}(x,y) \mid (x,y)\in \Cc^*\times \Cc^* \text{ and } 
\frac{\partial f_{\gamma_x}}{\partial x}(x,y) = \frac{\partial f_{\gamma_x}}{\partial y}(x,y) = 0 \bigg\rbrace.$$
In a similar way we define $\mathcal{C}_y$.

A result of \cite[Proposition 6]{NZ} is:
\begin{theorem}
\label{th:NZ}Let $f\in \Cc[x,y]$ be a non-degenerate and non-convenient polynomial with $f(0,0)=0$,
not depending only on one variable.
The set of critical values at infinity of $f$ is
$$\Binf=\mathcal{C}_x\cup \mathcal{C}_y \quad \text{ or } \quad \Binf = \{0\} \cup \mathcal{C}_x\cup \mathcal{C}_y.$$
\end{theorem}

Unfortunately this theorem does not determine whether $0 \in \Binf$
(and notice that the value $0$ may be already included in $\mathcal{C}_x$ or $\mathcal{C}_y$).
This value $0$ is treated in the following lemma. 
\begin{lemma}
\label{lem:0binf}

Let $f\in \Cc[x,y]$ be a non-degenerate and non-convenient polynomial, 
with isolated singularities and with $f(0,0)=0$. Then
 $$\Binf = \Binfx \cup \Binfy$$
where we define:
\begin{enumerate}
  \item if $f$ is convenient for the $x$-axis then $\Binfx:=\varnothing$;
  \item otherwise there exists $x^py$ in $\supp(f)$ where $p\ge 0$ is supposed to be maximal;
  \begin{enumerate}
    \item If $x^py$ is in a face of $\Gamma_-(f)$ then $\Binfx := \mathcal{C}_x$
and $0\notin \Binfx$;
    \item If $x^py$ is not in a face of $\Gamma_-(f)$ then $\Binfx := \{0\} \cup 
\mathcal{C}_x$;
  \end{enumerate}
  \item we set a similar definition for $\Binfy$.
\end{enumerate}
\end{lemma}

Theorem \ref{th:NZ} and its refinement Lemma \ref{lem:0binf} enable to calculate $\Binf$ from
$\supp(f)$. The different cases of Lemma \ref{lem:0binf} are pictured in Figures 
\ref{fig:gc1} and \ref{fig:gc2}.

\begin{proof}
As $f$ is non-convenient with $f(0,0)=0$ we may suppose that
$f$ is non-convenient for the $x$-axis so that 
$f(x,y)=yk(x,y)$. But $f$ has isolated 
singularities, so $y$ does not divide $k$. Then there is a monomial $x^py \in \supp(f)$,
we can suppose that $p\ge 0$ is maximal among monomials $x^ky \in \supp(f)$.

Let $d= \deg f$. Let $\bar f(x,y,z) -cz^d$ be the homogeneization of $f(x,y)-c$;
at the point at infinity $P = (1:0:0)$, we define $g_c(y,z) = \bar f(1,y,z) -cz^d$. 
Notice that only $(1:0:0)$ and $(0:1:0)$ can be singularities at infinity for $f$.
The value $0$ is a critical value at infinity for the point at infinity $P$ (that is to say $0\in \Binfx$) 
if and only if $\mu_{P}(g_0) > \mu_{P}(g_c)$
where $c$ is a generic value.

The Newton polygon of the germ of singularity $g_c$ can be computed from the Newton polygon
$\Gamma(f)$, for $c \not=0$, see \cite[Lemma 7]{NZ}.
If $A,B,O$ are the points on the Newton diagram of coordinates $(d,0), (0,d), (0,0)$,
then the Newton diagram of $g_c$ has origin $A$ with $y$-axis $AB$, $z$-axis $AO$, and
the convex closure of $\supp(g_c)$ corresponds to $\Gamma_-(f)$.

\begin{figure}[ht]
\begin{center}
\unitlength 0.8mm
\begin{picture}(70,40)(0,-3)

\put(0,0){\vector(1,0){70}}
\put(0,0){\vector(0,1){30}}

\put(-1,-2){\makebox(0,0){$O$}}

\put(0,20){\circle*{2}}
\put(0,20){\line(3,1){30}}
\put(30,30){\circle*{2}}
\put(30,30){\line(1,-1){10}}
\put(40,20){\circle*{2}}
\put(40,20){\line(-1,-1){10}}
\put(30,10){\circle*{2}}
\put(30,10){\line(-3,-1){30}}

\put(60,0){\vector(-1,1){40}}
\put(61,2){\makebox(0,0){$A$}}
\put(0,0){\vector(0,1){30}}
\put(20,36){\makebox(0,0){$y$}}
\put(60,0){\vector(-1,0){70}}
\put(-10,2){\makebox(0,0){$z$}}
\put(25,11){\makebox(0,0){$x^py$}}

\end{picture}
\caption{Newton polygon of $g_c$. First case: $0\not\in \Binfx$. \label{fig:gc1}}
\end{center}
\end{figure}

\begin{figure}[ht]
\begin{center}
\unitlength 0.8mm
\begin{picture}(70,36)(0,-3)

\put(0,0){\vector(1,0){70}}
\put(0,0){\vector(0,1){30}}
\put(-1,-2){\makebox(0,0){$O$}}

\put(0,20){\circle*{2}}
\put(0,20){\line(3,1){30}}
\put(30,30){\circle*{2}}
\put(30,30){\line(1,-1){10}}
\put(40,20){\circle*{2}}
\put(40,20){\line(-2,-1){40}}

\put(10,10){\circle*{2}}
\put(5,11){\makebox(0,0){$x^py$}}

\put(60,0){\vector(-1,1){40}}
\put(61,2){\makebox(0,0){$A$}}
\put(0,0){\vector(0,1){30}}
\put(20,36){\makebox(0,0){$y$}}
\put(60,0){\vector(-1,0){70}}
\put(-10,2){\makebox(0,0){$z$}}

\end{picture}

\caption{Newton polygon of $g_c$. Second case: $0\in \Binfx$. \label{fig:gc2}}
\end{center}
\end{figure}

We denote by $\Delta_c$ the Newton polygon of the germ $g_c$, for
a generic value $c$, $\Delta_{c}$ is non-degenerate and 
$\mu_P(g_c) = \nu(\Delta_{c})$.
The Newton polygon $\Delta_0$ has no common point with the $z$-axis $AO$
but $\nu$ may be defined for non-convenient series, see \cite[Definition 1.9]{K}.

If $x^py$ is in the face $\gamma_x$ of $\Gamma_-(f)$ then $\Delta_0$ is non-degenerate
and $ \nu(\Delta_0)=\nu(\Delta_c)$,
then by \cite[Theorem 1.10]{K}
 $\mu_P(g_0) = \nu(\Delta_0)$ and $\mu_P(g_c) = \nu(\Delta_c)$. 
So $\mu_P(g_0)=\mu_P(g_c)$ and $0$ is not a critical value at infinity for the point $P$ :
$0\notin \Binfx$.

If $x^py$ is not in a face of $\Gamma_-(f)$ then there is a triangle $\Delta_c$
that disappears in $\Delta_0$, by the positivity of $\nu$ (see below) we have
$\nu(\Delta_0) > \nu(\Delta_c)$,
then by \cite[Theorem 1.10]{K}:
$\mu_P(g_0) \ge \nu(\Delta_0) > \nu(\Delta_c) = \mu_P(g_c)$.
So we have $0 \in \Binfx$.
\qed
\end{proof}

\subsection{Additivity and positivity}

We need a variation of Kouchnirenko's number $\nu$.
Let $T$ be a polytope whose vertices are in $\Nn\times \Nn$,  $S>0$ the area of $T$, 
 $a$ the length of the intersection of $T$
with the $x$-axis, and $b$ the length of the intersection of $T$ with the $y$-axis.
 We define 
$$\tau(T) = 2S-a-b, \text { so that, } \nu(T) = \tau(T)+1.$$
It is clear that $\tau$ is additive:
$\tau(T_1 \cup T_2) = \tau(T_1) + \tau(T_2) - \tau(T_1\cap T_2)$, and
in particular if $T_1\cap T_2$ has null area then
$\tau(T_1 \cup T_2) = \tau(T_1) + \tau(T_2)$.
This formula enables us to argue on triangles only (after a triangulation of $T$).

Let $T_0$ be the triangle defined by the vertices $(0,0), (1,0), (0,1)$, we have $\nu(T_0)=-1$.
We have the following facts, for every triangle $T \not= T_0$:
\begin{enumerate}
  \item $\nu(T) \ge 0$;
  \item $\nu(T) = 0$ if and only if $T$ has an edge contained in the $x$-axis or the $y$-axis
and the height of $T$ (with respect to this edge) is $1$.
\end{enumerate}

\begin{remark*}
The formula of additivity can be generalized in the $n$-di\-men\-sion\-al case, but the positivity
can not. Here is a counter-example found by G\"{u}nter Ewald:
Let $n=4$,  $a$ a positive integer and let $T$ be the polytope whose vertices are:
$(1,0,0,0)$, $(1+a,0,0,0)$, $(1,1,1,0)$, $(1,2,1,0)$, $(1,1,1,1)$ then
$\tau(T)= \nu(T)+1 = -a < 0$. 
\end{remark*}

\subsection{Families of polytopes}

We consider a family $(f_s)_{s\in[0,1]}$ of complex polynomials in two variables 
with isolated singularities. We suppose that $\mu(s)+\lambda(s)$ remains constant.
We denote by $\Gamma(s)$ the Newton polygon of $f_s$.
We suppose that $f_s$ is non-degenerate for $s\in ]0,1]$.

\begin{figure}[ht]
\begin{center}
\unitlength 0.8mm
\begin{picture}(60,43)(0,-3)

\put(0,0){\vector(1,0){60}}
\put(0,0){\vector(0,1){40}}
\put(60,2){\makebox(0,0){$x$}}
\put(2,40){\makebox(0,0){$y$}}

\put(0,30){\circle{2}}
\put(0,30){\line(1,0){30}}
\put(0,10){\circle*{2}}
\put(0,10){\line(3,2){30}}
\put(30,30){\circle*{2}}
\put(30,30){\line(1,-1){20}}
\put(50,10){\circle{2}}
\put(50,10){\line(-2,-1){20}}
\put(30,0){\circle*{2}}
\put(30,0){\line(0,1){30}}

\put(35,15){\makebox(0,0){$T_2$}}
\put(6,22){\makebox(0,0){$T_1$}}

\end{picture}
\caption{Triangles that disappear.\label{fig:dis}}
\end{center}
\end{figure}
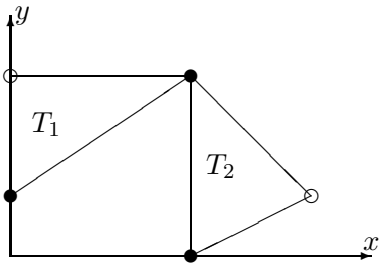

We will always assume that the only critical parameter is $s=0$. We will say that a monomial
$x^py^q$ \defi{disappears} if $(p,q) \in \supp(f_s) \setminus \supp(f_0)$ for $s\not=0$.
By extension a triangle of $\Nn\times \Nn$ disappears if one of its vertices 
(which is a vertex of $\Gamma(s)$, $s\not=0$) disappears.
Now after a triangulation of $\Gamma(s)$ we have a finite number of 
triangles $T$ that disappear (see Figure \ref{fig:dis}, on  pictures of the Newton diagram, 
a plain circle is drawn for a monomial that does not disappear and 
an empty circle for monomials that disappear).

We will widely use the following result, under the hypotheses of Theorem \ref{th:deg}:
\begin{lemma}
\label{lem:dis}
Let $T\not=T_0$ be a triangle that disappears then $\tau(T)=0$.
\end{lemma}

\begin{proof}
We suppose that $\tau(T)>0$.
By the additivity and positivity of $\tau$ we have for $s \in ]0,1]$:
$$\nu(s) = \nu(\Gamma(s)) \ge \nu(\Gamma(0))+\tau(T) > \nu(0).$$
Then by Theorem \ref{th:CN},
$$\mu(s)+\lambda(s) = \nu(s) > \nu(0) \ge \mu(0)+\lambda(0).$$
This gives a contradiction with $\mu(s)+\lambda(s)=\mu(0)+\lambda(0)$.

We remark that we do not need $f_0$ to be non-degenerate because in all cases we have
$\nu(0) \ge \mu(0)+\lambda(0)$.
\qed
\end{proof}

\section{Constancy of the degree}
\label{sec:demdeg}

\subsection{Families of quasi-constant degree}

Let $\sigma \in [0,1]$, we choose a small enough neighborhood $I$ of $\sigma$.
Let $\mathcal{M}_\sigma$ be the set of monomials that disappear at $\sigma$:
$\mathcal{M}_\sigma = \supp(f_s)\setminus \supp(f_\sigma)$ for $s\in I\setminus\{\sigma\}$.
The family $(f_s)_{s\in[0,1]}$ is of \defi{quasi-constant degree at $\sigma$}
if 
\begin{align*}
 \text{there exists } x^py^q &\in \supp(f_\sigma) \text{ such that } \\
 &(\forall x^{p'}y^{q'} \in \mathcal{M}_\sigma \ \ (p>p') \text{ or } (p=p' \text{ and } q>q')) \\
\text{ or } &(\forall x^{p'}y^{q'} \in \mathcal{M}_\sigma \ \ (q>q') \text{ or } (q=q' \text{ and } p>p')).
\end{align*}
The family $(f_s)_{s\in[0,1]}$ is of \defi{quasi-constant degree} if it is of quasi-constant degree
at each point $\sigma$ of $[0,1]$.
 The terminology is justified by the following remark:
\begin{lemma}
If $(f_s)$ is of quasi-constant degree at $\sigma\in [0,1]$, then there exists $\Phi \in \Aut \Cc^2$ such 
that $\deg  f_s \circ \Phi$ is constant in a neighborhood of $\sigma$.
\end{lemma}

The proof is simple: suppose that $x^py^q$ is a monomial of $\supp(f_\sigma)$
such that for all $ x^{p'}y^{q'} \in \mathcal{M}_\sigma$, $p>p'$ or $(p=p' \text{ and } q>q')$.
We set $\Phi(x,y) = (x+y^\ell,y)$ with $\ell \gg 1$. 
Then the monomial of highest degree in $f_s \circ \Phi$ is $y^{q+p\ell}$ and 
does not disappear at $\sigma$.
For example let $f_s(x,y) = xy+sy^3$, we set 
$\Phi(x,y) = (x+y^3,y)$ then $f_s \circ \Phi (x,y) = y^4+xy+sy^3$ is of constant degree.

\bigskip

We prove Theorem \ref{th:deg}. 
We suppose that the degree changes, more precisely we suppose that
$\deg f_s$ is constant for $s\in ]0,1]$ and that $\deg f_0 < \deg f_s$,
$s\in ]0,1]$. As the degree changes the Newton polygon $\Gamma(s)$ cannot be constant, that means that
at least one vertex of $\Gamma(s)$ disappears.

\subsection{Exceptional case}
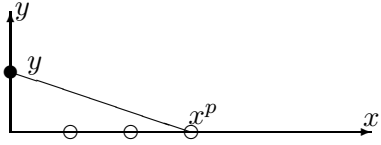
\begin{figure}[ht]
\begin{center}
\unitlength 0.8mm
\begin{picture}(60,23)(0,-3)

\put(0,0){\vector(1,0){60}}
\put(0,0){\vector(0,1){20}}
\put(60,2){\makebox(0,0){$x$}}
\put(2,20){\makebox(0,0){$y$}}

\put(0,10){\circle*{2}}
\put(30,0){\circle{2}}
\put(20,0){\circle{2}}
\put(10,0){\circle{2}}
\put(0,10){\line(3,-1){30}}

\put(4,11){\makebox(0,0){$y$}}
\put(32,3){\makebox(0,0){$x^p$}}

\end{picture}

\caption{Case $f_0 \in \Cc[y]$. \label{fig:exc}}
\end{center}
\end{figure}

We suppose that $f_0$ is a one-variable polynomial, for example $f_0 \in \Cc[y]$.
As $f_0$ has isolated singularities then $f_0(x,y) = a_0y+b_0$, so $\mu(0)=\lambda(0)=0$,
then for all $s$, $\mu(s)=\lambda(s)=0$. So $\nu(s) = \nu(\Gamma(s))=0$, then
$\deg_y f_s = 1$, and $f_s(x,y)= a_sy+b_s(x)$, so $(f_s)_{s\in[0,1]}$ is 
a family of quasi-constant degree (see Figure \ref{fig:exc}).
We exclude this case for the end of the proof.

\subsection{Case to exclude}
\begin{figure}[ht]
\begin{center}
\unitlength 0.8mm
\begin{picture}(60,38)(0,-3)

\put(0,0){\vector(1,0){60}}
\put(0,0){\vector(0,1){35}}
\put(60,2){\makebox(0,0){$x$}}
\put(2,35){\makebox(0,0){$y$}}

\put(0,10){\circle*{2}}
\put(0,10){\line(1,1){20}}
\put(20,30){\circle*{2}}
\put(20,30){\line(1,-1){10}}
\put(30,20){\circle{2}}
\put(30,20){\line(-1,-2){10}}
\put(20,0){\circle*{2}}

\put(20,30){\line(0,-1){30}}

\put(24,17){\makebox(0,0){$T$}}
\put(36,20){\makebox(0,0){$x^py^q$}}
\end{picture}
\caption{Case where a monomial $x^py^q$, $p>0,q>0$ of $\Gamma(s)$ disappears. \label{fig:dis2}}
\end{center}
\end{figure}
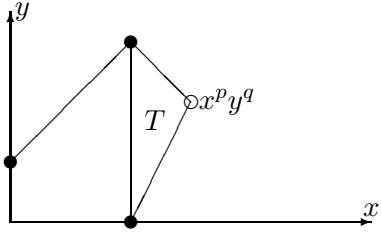

We suppose that a vertex $x^py^q$, $p>0,q>0$ of $\Gamma(s)$ disappears.
Then there exists a triangle $T$ that disappears whose faces are not contained
in the axis. Then $\tau(T)>0$ that contradicts Lemma \ref{lem:dis} (see Figure \ref{fig:dis2}).

\subsection{Case where a monomial $x^a$ or $y^b$ disappears (but not both)}

If, for example the monomial $y^b$ of $\Gamma(s)$ disappears and $x^a$ does not,
then we choose a monomial $x^py^q$, with maximal $p$, among monomials in $\supp(f_s)$.
Certainly $p\ge a>0$. We also suppose that $q$ is maximal among monomials $x^py^k \in \supp(f_s)$.
If $q=0$ then $p=a$, and the monomial $x^py^q=x^a$ does not disappear (by assumption).
If $q>0$ then $x^py^q$ cannot disappear (see above). In both cases the monomial $x^py^q$ proves
that  $(f_s)$ is of quasi-constant degree.

\subsection{Case where both $x^a$ and $y^b$ disappear}

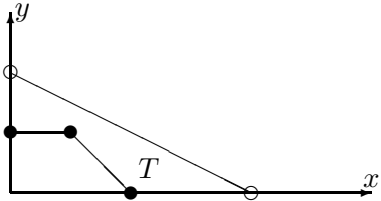
\begin{figure}[ht]
\begin{center}
\unitlength 0.8mm
\begin{picture}(60,33)(0,-3)

\put(0,0){\vector(1,0){60}}
\put(0,0){\vector(0,1){30}}
\put(60,2){\makebox(0,0){$x$}}
\put(2,30){\makebox(0,0){$y$}}

\put(0,20){\circle{2}}
\put(0,20){\line(2,-1){40}}
\put(40,0){\circle{2}}

\put(0,10){\circle*{2}}
\put(0,10){\line(1,0){10}}
\put(10,10){\circle*{2}}
\put(10,10){\line(1,-1){10}}
\put(20,0){\circle*{2}}

\put(23,4){\makebox(0,0){$T$}}

\end{picture}

\caption{Sub-case : no monomial $x^py^q$ in $\Gamma(s)$, $p>0,q>0$. \label{fig:subcase1}}
\end{center}
\end{figure}

\paragraph{Sub-case : No monomial $x^py^q$ in $\Gamma(s)$, $p>0, q>0$.} 
Then there is an area $T$ with $\tau(T) > 0$ 
that disappears (see Figure \ref{fig:subcase1}). Contradiction.

\begin{figure}[ht]
\begin{center}
\unitlength 0.8mm
\begin{picture}(60,53)(0,-3)

\put(0,0){\vector(1,0){40}}
\put(0,0){\vector(0,1){50}}
\put(42,2){\makebox(0,0){$x$}}
\put(2,50){\makebox(0,0){$y$}}

\put(0,40){\circle{2}}
\put(0,40){\line(1,-1){10}}
\put(10,30){\circle*{2}}
\put(10,30){\line(1,-2){10}}
\put(20,10){\circle*{2}}
\put(20,10){\line(0,-1){10}}
\put(20,0){\circle{2}}

\put(0,20){\circle*{2}}
\put(10,30){\line(-1,-1){10}}

\put(10,0){\circle*{2}}
\put(20,10){\line(-1,-1){10}}

\put(17,3){\makebox(0,0){$T_2$}}
\put(4,29){\makebox(0,0){$T_1$}}

\put(4,40){\makebox(0,0){$y^b$}}
\put(15,30){\makebox(0,0){$xy^{q'}$}}
\put(25,10){\makebox(0,0){$x^py$}}
\put(24,2){\makebox(0,0){$x^a$}}

\end{picture}
\caption{Sub-case : existence of monomials $x^py^q$ in $\Gamma(s)$, $p>0,q>0$.\label{fig:slope}}
\end{center}
\end{figure}
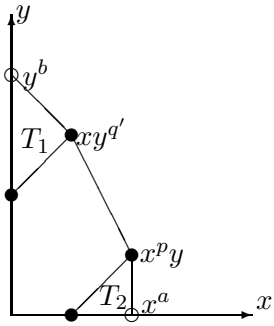

\paragraph{Sub-case : there exists a monomial $x^py^q$ in $\Gamma(s)$, $p>0,q>0$.}
We know that $x^py^q$ is in $\Gamma(0)$ because it cannot disappear.
As $\deg f_0 < \deg f_s$, a monomial $x^py^q$ that does not disappear verifies
$\deg x^py^q = p+q < \deg f_s$, ($s\in]0,1]$).
So the monomial of highest degree is $x^a$ or $y^b$. We will suppose that it is $y^b$,
so $d=b$, and the monomial $y^b$ disappears. Let $x^py^q$ be a monomial of $\Gamma(s)$, $p,q>0$
with minimal $q$. By assumption such a monomial exists.
Then certainly we have $q=1$, otherwise there exists a region $T$ that disappears
 with $\tau(T)>0$ (on Figure \ref{fig:slope} the regions $T_1$ and $T_2$ verify $\nu(T_1)=0$
and $\nu(T_2)=0$).
For the same reason the monomial $x^{p'}y^{q'}$ with minimal $p'$ verifies $p'=1$.

We look at the segments of $\Gamma(s)$, starting from $y^b=y^d$ and ending
at $x^a$. The first segment is from $y^d$ to $xy^{q'}$, $(p'=1)$ and we know that $p'+q'<d$ so 
the slope of this segment is strictly less than $-1$. By the convexity of $\Gamma(s)$ all the following
slopes are strictly less than $-1$.
The last segment is from $x^py$ to $x^a$, with a slope strictly less than  $-1$, so $a \le p$.
Then the monomial $x^py$ gives that  $(f_s)_{s\in[0,1]}$ is of quasi-constant degree.

\section{Continuity of the critical values}

We now prove Theorem \ref{th:cont}. 
We will suppose that $s=0$ is the only problematic parameter. In particular
$\Gamma(s)$ is constant for all $s\in ]0,1]$.

\subsection{The Newton polygon changes} 
That is to say $\Gamma(0) \not= \Gamma(s)$, $s\not=0$.
As in the proof of Theorem \ref{th:deg} (see paragraph \ref{sec:demdeg}) we remark:
\begin{itemize}
  \item If $f_0$ is a one-variable polynomial then $\Binf(s)=\varnothing$ for all $s\in [0,1]$.
  \item A vertex $x^py^q$, $p>0, q>0$ of $\Gamma(s)$ cannot disappear.
\end{itemize}

So we suppose that a monomial $x^a$ of $\Gamma(s)$ disappears
(a similar proof holds for $y^b$). 
Then for $s \in ]0,1]$ the monomial $x^a$ is in $\Gamma(s)$, so
there are no critical values at infinity for $f_s$ at the point $P=(1:0:0)$.
If $\Gamma(0)$ contains a monomial $x^{a'}$, $a'>0$ then there are no critical values
 at infinity for $f_0$ at the point $P$ .
So we suppose that all monomials $x^k$ disappear. Then a monomial $x^py^q$ of $\supp(f_0)$ with 
minimal $q>0$, verifies $q=1$, otherwise there would exist a region $T$ with $\tau(T)>0$
(in contradiction with the constancy of $\mu(s)+\lambda(s)$, see Lemma \ref{lem:dis}).
And for the same reason if we choose $x^py$ in $\supp(f_0)$ with maximal $p$ then $p>0$ and 
 $x^py \in \Gamma(0)$.
Now the edge of $\Gamma_-(f_0)$ that contains the origin and the monomial $x^py$ (with maximal $p$)
begins at the origins and ends at $x^py$ (so in particular there is no monomial $x^{2p}y^2$, 
$x^{3p}y^3$ in $\supp(f_0)$). Now from Theorem \ref{th:NZ} and Lemma \ref{lem:0binf}
we get that there are no critical values at infinity for $f_0$ at $P$.

So in case where $\Gamma(s)$ changes, we have for all $s\in[0,1]$, $\Binf(s) = \varnothing$.

\subsection{The Newton polygon is constant : case of non-zero critical values} 

We now prove the following lemma that ends the proof of Theorem \ref{th:cont}.
\begin{lemma}
Let a family $(f_s)_{s\in[0,1]}$ such that
$f_s$ is non-degenerate for all $s\in[0,1]$ and $\Gamma(s)$ is constant, then
the multi-valued function $s\mapsto \Binf(s)$ is continuous and closed.
\end{lemma}

In this paragraph and the next one we suppose that $f_s(0,0)=0$, 
that is to say the constant term of $f_s$ is zero.
We suppose that $c(0) \in \Binf(0)$ and that $c(0) \not=0$.
Then $c(0)$ has been obtained by the result of N\'em\'ethi-Zaharia
(see Theorem \ref{th:NZ}).
There is a face $\gamma$ of $\Gamma_-(f_0)$ that contains the origin
such that $c(0)$ is in the set: 
$$\mathcal{C}_\gamma(0) = \bigg\lbrace (f_0)_\gamma(x,y) \mid (x,y)\in (\Cc^*)^2 \text{ and } 
\frac{\partial (f_0)_\gamma}{\partial x}(x,y) = \frac{\partial (f_0)_\gamma}{\partial y}(x,y) = 0 \bigg\rbrace.$$
Now, as $\Gamma(s)$ is constant, $\gamma$ is a face of $\Gamma_-(s)$ for all $s$.
There exists a family of polynomials $h_s \in \Cc[t]$ and a monomial $x^py^q$ ($p,q>0$, $\gcd(p,q)=1$) such that
$(f_s)_\gamma(x,y) = h_s(x^py^q)$. The family $(h_s)$ is continuous (in $s$) and is 
of constant degree (because $\Gamma(s)$ is constant).
The set $\mathcal{C}_\gamma(0)$ and more generally the set  $\mathcal{C}_\gamma(s)$
can  be computed by
$$\mathcal{C}_\gamma(s) = \bigg\lbrace h_s(t) \mid  t \in \Cc^* \text{ and }
h_s'(t)=0  \bigg\rbrace.$$
As $c(0) \in \mathcal{C}_\gamma(0)$ there exists a $t_0\in \Cc^*$ with $h_0'(t_0)=0$,
and for $s$ near $0$ there is a $t_s\in \Cc^*$ near $t_0$ with $h_s'(t_s)=0$
(because $h'_s(t)$ is a continuous function of $s$ of constant degree in $t$).
Then $c(s) = h_s(t_s)$ is a critical value at infinity near $c(0)$
and we get the continuity.

\subsection{The Newton polygon is constant : case of the value $0$}

We suppose that $c(0)=0 \in \Binf(0)$ and that $f_s(x,y)=yk_s(x,y)$.
We will deal with the point at infinity $P=(1:0:0)$, the point
$(0:1:0)$ is treated in a similar way.
Let $x^py$ be a monomial of $\supp(f_s)$ with maximal $p \ge 0$, $s\not=0$.
If $x^py$ is not in a face of $\Gamma(s)$ then $0 \in \Binf(s)$  for all $s\in [0,1]$,
and we get the continuity. Now we suppose that
$x^py$  is in a face of $\Gamma(s)$; then  $x^py$ disappears
otherwise $0$ is not a critical value
at infinity (at the point $P$) for all $s\in [0,1]$.
As $\Gamma(s)$ is constant then the face $\gamma$ that contains the origin and $x^py$ 
for $s\not=0$ is also a face of $\Gamma(0)$, then there exists a monomial 
$(x^py)^k$, $k>1$ in $\supp(f_0)$.
Then $(f_s)_\gamma = h_s(x^py)$, $h_s  \in \Cc[t]$. We have $\deg h_s > 1$, with
$h_s(0)=0$ (because $f(0,0)=0$) and $h_0'(0)=0$ (because $x^py$ disappears). Then
$0 \in \mathcal{C}_\gamma(0) \subset \Binf(0)$ but by continuity of $h_s$ we have a critical value
$c(s) \in \mathcal{C}_\gamma(s) \subset \Binf(s)$ such that $c(s)$ tends towards $0$ (as $s\rightarrow 0$).
It should be noticed that for $s\not=0$, $c(s)\not=0$.

In all cases we get the continuity of $\Binf(s)$.

\subsection{Proof of the closeness of $s\mapsto \Binf(s)$}

We suppose that $c(s) \in \Binf(s)$, is a continuous function of $s\not=0$, with
a limit $c(0)\in \Cc$ at $s=0$. We have to prove that $c(0) \in \Binf(0)$.
As there are critical values at infinity we suppose that $\Gamma(s)$ is constant.

\paragraph{Case $c(0)\not=0$.} Then for $s$ near $0$, $c(s) \not=0$ by continuity,
then $c(s)$ is obtained as a critical value of $h_s(t)$. By continuity
$c(0)$ is a critical value of $h_0(t)$: $h_0(t_0)=c(0)$, $h_0'(t_0)=0$;
as $c(0)\not=0$, $t_0\not=0$ (because $h_0(0)=0$). Then $c(0) \in \Binf(0)$.

\paragraph{Case $c(0)=0$.} Then  let $x^py$ be the monomial of $\supp(f_s)$, $s\not=0$, with
maximal $p$. By Lemma \ref{lem:0binf} if $x^py \notin \Gamma(s)$ for $s\in]0,1]$ then
$0\in \Binf(s)$ for all $s\in[0,1]$ and we get closeness. If $x^py \in \Gamma(s)$, $s\not=0$,
then as $c(s) \rightarrow 0$ we have that $x^py$ disappears, so $x^py \notin \Gamma(0)$,
then by Lemma \ref{lem:0binf}, $c(0)=0 \in \Binf(0)$.

\subsection{Proof of the closeness of $s\mapsto \B(s)$}

We now prove Corollary \ref{th:contb}.
The multi-valued function $s\mapsto \B(s)$ is continuous because $\B(s) = \Baff(s) \cup \Binf(s)$
and $s \mapsto \Baff(s)$, $s\mapsto \Binf(s)$ are continuous.
For closeness, it remains to prove that if $c(s) \in \Baff(s)$ is a continuous function
with a limit $c(0)\in \Cc$ at $s=0$ then $c(0) \in \B(0)$.

We suppose that $c(0) \notin \Baff(0)$. There exist critical points $Q_s = (x_s,y_s) \in \Cc^2$ 
of $f_{s}$ with $f_s(x_s,y_s) = c(s)$, $s \not=0$.
We can extract a countable set $\mathcal{S}$ of $]0,1]$ such that the sequence
$(Q_s)_{s\in \mathcal{S}}$ converges towards $P$ in $\Cc P^2$.
As $c(0) \notin \Baff(0)$ we have that $P$ relies on the line at infinity and
we may suppose that $P=(0:1:0)$.

By Theorem \ref{th:deg} we may suppose, after an algebraic automorphism of $\Cc^2$ if necessary,
that $d=\deg f_s$ is constant.
Now we look at $g_{s,c}(x,z) = \bar{f}_s(x,1,z)-cz^d$.
The critical point $Q_s$ of $f_s$ with critical value $c(s)$ gives
a critical point $Q'_s = (\frac{x_s}{y_s},\frac{1}{y_s})$ of $g_{s,c(s)}$ with critical value $0$
(see \cite[Lemma 21]{Bo}).
Then by semi-continuity of the local Milnor number on the fiber $g_{s,c(s)}^{-1}(0)$ we
have $\mu_P(g_{0,c(0)}) \ge \mu_P(g_{s,c(s)}) + \mu_{Q'_s}(g_{s,c(s)}) >  \mu_P(g_{s,c(s)})$.
As $\mu(s) +\lambda(s)$ is constant we have  $\mu_P(g_{s,c})$ constant for a generic $c$
(see \cite[Corollary 5.2]{ST} or \cite{BT}).
Then we have $\mu_P(g_{0,c(0)}) - \mu_P(g_{0,c}) > \mu_P(g_{s,c(s)})-\mu_P(g_{s,c}) \ge 0$.
Then $c(0) \in \Binf(0)$.
And we get closeness for $s \mapsto \B(s)$.


\end{document}